\input amstex
\input amsppt.sty

\magnification=1200
\voffset=-1.5cm

\NoBlackBoxes
\TagsOnRight
\NoRunningHeads

\def\d{\downharpoonright}
\def\deg{\operatorname{deg}}
\def\sgn{\operatorname{sgn}}
\def\bZ{\Bbb Z}
\def\l{\langle}
\def\r{\rangle}

\def\Res{\text{\rm Res}\,}

\def\Om{\Omega}
\def\la{\lambda}
\def\om{\omega}
\def\ga{\gamma}
\def\ep{\epsilon}

\def\si{\sigma}

%%%%%%%%%%%%%%%%%%%%%%%%%%%%%%%%%%%%%%%%%%%%%%%%%%%%%
%%%  Picture drawing routines by A.~Yu.~Okounkov  %%%
%%%%%%%%%%%%%%%%%%%%%%%%%%%%%%%%%%%%%%%%%%%%%%%%%%%%%
\def\li #1 #2 #3 {\vrule height #1  width #2  depth #3 }
% \li draws a black box

\def\hl #1 #2 #3 {\rlap{ \kern #1 pc \raise -#2 pc \hbox{ \li {0.1pt} {#3pc} 0pt }}}
% \hl draws a horizontal line starting in the point (#1,#2)
% with length #3 and width #4

\def\vl #1 #2 #3 {\rlap{ \kern #1 pc \raise -#2 pc \hbox{ \li {#3pc} {0.1pt} 0pt }}}
% \vl draws a vertical line starting in the point (#1,#2)
% with length #3 and width #4

\def\bhl #1 #2 #3 {\rlap{ \kern #1 pc \raise -#2 pc \hbox{ \li {1pt} {#3pc} 0pt }}}
% \hl draws a horizontal line starting in the point (#1,#2)
% with length #3 and width #4

\def\bvl #1 #2 #3 {\rlap{ \kern #1 pc \raise -#2 pc \hbox{ \li {#3pc} {1pt} 0pt }}}
% \vl draws a vertical line starting in the point (#1,#2)
% with length #3 and width #4

\def\wr #1 #2 #3 {\rlap{ \kern #1 pc \raise -#2 pc \hbox{\tenrm #3 }}}
% \wr writes #3 starting from the point (#1,#2)

\def\boxit #1 #2 #3 {\rlap{ \kern #1 pc \raise -#2 pc \vbox to 2 pc{
\vfill \hbox to 2pc{\hfill {\tenrm #3} \hfill } \vfill }}}

\def\botit #1 #2 #3 {\rlap{ \kern #1 pc \raise -#2 pc \hbox{\raise 4 pt
\hbox to 2pc{\hfill {\tenrm #3} \hfill }}}}

\def\picture #1 #2 #3 \endpicture
   {\smallskip \centerline{\hbox to #1 pc
   {\nullfont #3 \hss}} \nobreak \smallskip \centerline{#2} \smallskip}
% creates a picture with width #1 and name #2

\def\bx #1 #2 #3 {\boxit #1 #2 {$#3$} }

\topmatter

\title
{THE DIMENSION OF SKEW SHIFTED YOUNG DIAGRAMS, AND PROJECTIVE CHARACTERS OF
THE INFINITE SYMMETRIC GROUP}
\endtitle

\author
{Vladimir~N.~Ivanov}
\endauthor

\address
Russia, 119992, Moscow, Vorobjevy
gory, Moscow State University, GZ,
Department of Mathematics and
Mechanics, Chair of Higher Algebra.
\endaddress
\email
vivanov\@mccme.ru
\endemail
\endtopmatter

\document

\head \S0. Introduction
\endhead

This article
was originally published in Russian
in
"Representation Theory, Dynamical Systems, Combinatorial and
Algorithmic Methods. Part 2" (A.~M.~Vershik, ed.), Zapiski Nauchnyh
Seminarov POMI {\bf 240} (1997), 115--135 (this text in Russian is
available via
{\tt http://www.pdmi.ras.ru/znsl/1997/v240.html}).
As it was mentioned in the "Journal-Ref" field this English translation
was published in Journal of Mathematical Sciences (New York)
{\bf 96} (1999), no.~5, 3517--3530.

The dimension of a given skew shifted Young diagram is the number of
standard labellings of this diagram. In \S1 of this paper, we obtain a formula
for the dimension of an arbitrary skew shifted Young diagram. For this
purpose, we introduce polynomials
$P^*_\mu$
that are factorial analogues of a particular case of Hall--Littlewood
polynomials
$P(\cdot,t)$
for $t=-1$ (\cite{3, ch.~III, \S1}). The definition of the polynomials
$P^*_\mu$
is due to A.~Yu.~Okounkov.

As an application of the formula for the
dimension of a skew shifted Young
diagram, we obtain
new proof of the classification of projective characters of the group
$S(\infty)$. The classic Thoma's work \cite{12} contains the description
of characters (in von Neumann's sense) of the infinite symmetric
group $S(\infty)$ that is the inductive limit of the chain
of finite symmetric groups
$S(1)\subset S(2)\subset\ldots$
(characters in von Neumann's sense correspond
to finite factor--representations). In \cite{5} M.~L.~Nazarov
extended Thoma's theorem to projective
characters of the infinite symmetric
group $S(\infty)$. In \S2 of this paper, we give new proofs of
Nazarov's (\cite{5}) main results.

Earlier in \cite{6} and \cite{8} the analogous results in the
ordinary (non-projective) case were obtained. In this work
we follow the methods of \cite{8}.

The author is very grateful to G.~I.~Olshanski for setting the problem,
constant attention to the work and remarks on projects of the manuscript,
and to M.~L.~Nazarov for Remark~1.7 on the formula for the dimension
of a skew shifted diagram.

\head \S1.
Formula for the dimension of a skew shifted Young diagram
\endhead

A polynomial
$f(x_1,\ldots,x_n)$
is said to be {\it supersymmetric} if it satisfies the following conditions:
\roster
\item[1] $f$ is a symmetric polynomial in
$x_1,\ldots,x_n$;
\item[2] for any integers $i$, $j$ such that $1\le i<j\le n$, the polynomial
$$
f(x_1,\ldots,x_{i-1},t,x_{i+1},\ldots,x_{j-1},-t,
x_{j+1},\ldots,x_n)
$$
does not depend on $t$.
\endroster

Supersymmetric polynomials in $n$ variables form an algebra. We denote it
by $\Om(n)$.
It is graded by degrees of polynomials, i.e.,
$$
\Om(n)=\underset{k\ge0}\to\oplus\Om^k(n),
$$
where $\Om^k(n)$
consists of homogeneous polynomials of degree $k$ (including the zero
polynomial).

For $m\ge n$, we consider a homomorphism
$$
p_{m,n}:\Om(m)\to\Om(n)
$$
such that
$$
(p_{m,n}f)(x_1,\ldots,x_n)=f(x_1,\ldots,x_n,0,0,\ldots,0).
$$
Restricting to
$\Om^k(m)$,
we obtain linear transformations
$$
p^k_{m,n}:\Om^k(m)\to\Om^k(n).
$$

The projective limit of graded algebras
$\Om(n)$
in the category of graded algebras, taken
with respect to the morphisms $p_{m,n}$,
is an algebra and it is called the {\it algebra of supersymmetric
functions}. We denote it by
$\Om$.
By definition,
an element $f$ of $\Om$
is
a sequence
$(f_n)_{n\ge1}$ satisfying
the following conditions:

1) $f_n\in\Om(n)$, $n=1,2,\ldots$,

2) $f_{n+1}(x_1,\ldots,x_n,0)=f_n(x_1,\ldots,x_n)$ (stability condition),

3) $\sup\limits_n\deg f_n<\infty$.

Now we consider examples of supersymmetric functions that are important
for the sequel.

Let
$$
p_{k|n}(x_1,\ldots,x_n)=x^k_1+\ldots+x^k_n\,\,(k=1,2,\ldots).
$$
If
$\la=(\la_1,\ldots,\la_l)$ is a partition, then
$p_{\la|n}$ is defined as
$$
p_{\la|n}=p_{\la_1|n}\cdot p_{\la_2|n}\cdot\ldots
\cdot p_{\la_l|n}.
$$
If $\la$ is a partition such that all its
nonzero parts are odd, then
$$
p_{\la|n}\in\Om(n),
$$
and the sequence
$(p_{\la|n})_{n\ge1}$ defines a supersymmetric function $p_{\la}$.
Such functions $p_{\la}$
form a linear basis of the algebra
$\Om$. In other words,
$p_1,p_3,p_5,\ldots$, which are called odd Newton sums,
generate algebraically the algebra
$\Om$, see \cite{9},\cite{3}.

A partition is called {\it strict} if all its
nonzero parts are distinct.
The set of strict partitions of $n$ is denoted by $DP_n$.

Let
$\la=(\la_1,\la_2,\ldots)$ be an arbitrary partition. We denote by
$l(\la)$
the {\it length} of $\la$, i.e.,
the number of its nonzero parts. We denote by
$|\la|$ the {\it weight} of $\la$,
$$
|\la|=\la_1+\ldots+\la_{l(\la)}.
$$
We also use the notation
$\la\vdash n$ if $n=|\la|$.

In what follows, $\mu$ and $\la$ denote strict partitions, unless otherwise
specified.

First, we prove a proposition whose results we use in the sequel.

\proclaim{Proposition 1.1}{\it Let $x_1,x_2,\ldots$ be variables,
$r(x_1,\ldots,x_l)$ be a polynomial in $l$ variables. For
$n\ge l$, let
$$
R_n(x_1,\ldots,x_n)=r(x_1,\ldots,x_l)
\prod\limits_{i\le l,i<j\le n}
\frac{x_i+x_j}{x_i-x_j}
$$
and
$$
\widetilde R_n(x_1,\ldots,x_n)=\sum\limits_{\om\in S(n)}
R_n(x_{\om(1)},\ldots ,x_{\om(n)}).
$$
Then

{\rm a)} $\widetilde R_n$ is a polynomial, and
$$
\deg\widetilde R_n\le\deg r;
$$

{\rm b)} $\widetilde R_n$ is supersymmetric;

{\rm c)} $\widetilde R_n=0$, if $r$
is symmetric with respect to at least two variables
$x_i,x_j$, $i<j\le l$;

{\rm d)} if $x_1\ldots x_l$ divides $r(x_1,\ldots,x_l)$, then
$$
\widetilde R_{n+1}(x_1,\ldots ,x_n,0)=(n+1-l)\widetilde R_n(x_1,\ldots,
x_n).
$$
}\endproclaim

\demo{Proof} a) We will represent $\widetilde R_n$ as a ratio of
two polynomials. Denote by \linebreak
$V(x_1,\ldots,x_n)$ the Vandermonde determinant
$$
V(x_1,\ldots,x_n)=\prod\limits_{i<j}(x_i-x_j).
$$
We also set
$$
u_n(x_1,\ldots,x_n)=r(x_1,\ldots,x_l)
\prod\limits_{i\le l,i<j\le n}
(x_i+x_j)
\prod\limits_{l<i<j\le n}(x_i-x_j)
$$
and
$$
\widetilde u_n(x_1,\ldots,x_n)=
\sum\limits_{\om\in S(n)}\sgn(\om) u_n(x_{\om(1)},\ldots,
x_{\om(n)}).
\tag1.1
$$
Note the inequality
$$
\deg\widetilde u_n\le \deg V+\deg r.
\tag1.2
$$
We have the following relation between the polynomials
$$
\widetilde R_n=\frac{\widetilde u_n}V.
$$
It follows from (1.1) that $\widetilde u_n$ is a skew--symmetric polynomial
in $x_1,\ldots,x_n$,
hence $\widetilde R_n$ is a polynomial in $x_1,\ldots,x_n$.
It follows from (1.2) that
$$
\deg\widetilde R_n\le\deg r.
$$

b) Symmetry of $\widetilde R_n$
follows from its definition.

Let $x_i=t$, $x_j=-t$, for arbitrary integers
$i$ and $j$ such that $1\le i<j\le n$. Then
$R_n(x_{\om(1)},\ldots,x_{\om(n)})$ does not depend on
$t$ for any permutation $\om$ from the group $S(n)$. Hence, the polynomial
$\widetilde R_n$ also does not depend on $t$.

Thus, the polynomial $\widetilde R_n$ is supersymmetric.

c) Let $r$ be symmetric with respect to variables
$x_i$ and $x_j$, $i<j\le l$. Then the above--defined polynomial
$u_n$ is also symmetric with respect to $x_i$ and $x_j$. It follows
from (1.1) that
in this case
$\widetilde u_n=0$, since the sum in
(1.1) can be broken down into pairs of summands with equal
absolute values which occur
in this sum with different signs. Hence,
$\widetilde R_n=0$.

d) Let $x_1\ldots x_l$
divide $r(x_1,\ldots,x_l)$, and $x_{n+1}=0$.
Consider an arbitrary permutation $\om$ from the group
$S(n+1)$. If $\om^{-1}(n+1)\le l$, then
$$
R_{n+1}(x_{\om(1)},\ldots,x_{\om(n)}, x_{\om(n+1)})=0.
$$
If $\om^{-1}(n+1)>l$, then
$$
R_{n+1}(x_{\om(1)},\ldots,x_{\om(n)}, x_{\om(n+1)})=
R_n(x_{\si(1)},\ldots,x_{\si(n)}),
$$
where the permutation $\si$ from the subgroup $S(n)\subset S(n+1)$
is obtained from $\om$ by multiplying by the transposition
$(n+1,\om(n+1))$ from the left side.

This implies the desired equality
$$
\widetilde R_{n+1}=(n+1-l)\widetilde R_n.
$$
Proposition is proved. $\square$
\enddemo

Let us consider two particular cases of Proposition~1.1.

Let $\la$ be a partition, $l(\la)=l<n$,
$$
r(x_1,\ldots,x_l)=\frac{\prod\limits^l_{i=1}x^{\la_i}_i}
{(n-l)!}.
$$
In this case, the polynomial $\widetilde R_n$, which occurs in
the formulation of Proposition~1.1,
is denoted by $P_{\la|n}$. This is a particular case of
Hall--Littlewood polynomial when parameter
$t$ equals $-1$ ([3, ch.~III, \S1]). Put
$P_{\la|n}=0$ if $l(\la)>n$. It follows from Proposition~1.1 that
the sequence $(P_{\la|n})_{n\ge1}$ defines a supersymmetric function
$P_\la$. The functions $P_\la$
($\la$ denotes a strict partition),
which are called {\it Schur $P$-functions}, form a linear basis of
the algebra $\Om$, see [9]. Note that if
$\nu$ is a non--strict partition, then Proposition~1.1~c)
implies $P_\nu=0$.

We define the $k$th {\it decreasing factorial power} of a variable
$x$ as
$$
(x\downharpoonright k)
=\prod\limits^k_{i=1}(x-i+1),\quad k=1,2,\ldots.
$$
We also assume
$$
(x\d0)=1.
$$

Now we introduce the polynomials that play an important
role in \S1.

\proclaim{Definition 1.2 {\rm (A.~Yu.~Okounkov)}}{\it Let
$l=l(\la)\le n$,
$x_1,\ldots,x_n$ be variables. Let
$$
F_{\la|n}(x_1,\ldots,x_n)=\prod\limits^l_{i=1}
(x_i\d\la_i)
\prod\limits_{i\le l,i<j\le n}
\frac{x_i+x_j}{x_i-x_j}.
$$
We introduce a polynomial $P^*_{\la|n}$ by the formula
$$
P^*_{\la|n}=\frac1{(n-l)!}\sum\limits_{\om\in S(n)}
F_{\la|n}(x_{\om(1)},\ldots,x_{\om(n)}).
$$
}\endproclaim

\proclaim{Proposition 1.3}{\it If $l(\la)\le n$, then
$$
P^*_{\la|n}(x_1,\ldots,x_n)=P_{\la|n}(x_1,\ldots,x_n)+
g(x_1,\ldots,x_n),
$$
where $g(x_1,\ldots,x_n)$ is a supersymmetric polynomial of degree
less than $|\la|$.
}\endproclaim

\demo{Proof} If we set, in Proposition~1.1,
$$
r(x_1,\ldots,x_{l(\la)})=
\prod\limits^{l(\la)}_{i=1}(x_i\d\la_i)-
\prod\limits^{l(\la)}_{i=1}x^{\la_i}_i,
$$
then the polynomial $\widetilde R_n$ obtained in this Proposition
coincides (up to a scalar factor) with the difference
$$
P^*_{\la|n}-P_{\la|n}=g.
$$
By Proposition~1.1~a), we have
$$
\deg g\le\deg r<|\la|.\qquad \square
$$
\enddemo

If $l(\la)>n$, then put $P^*_{\la|n}=0$. It follows from
Proposition~1.1d) that the sequence
$(P^*_{\la|n})_{n\ge1}$ defines the supersymmetric function
$P^*_\la$.

Proposition~1.3 yields the form of the highest term of
$P^*_\la$.

\proclaim{Corollary 1.4}{\it
$$
P^*_\la=P_\la+g,
$$
where $g$ is a supersymmetric function of degree less than
$|\la|$.
}\endproclaim

Let $\mu$ be a strict partition. Denote
$$
H(\mu)=\prod\limits^{l(\mu)}_{t=1}
\mu_t!\prod\limits_{i<j}\frac{\mu_i+\mu_j}
{\mu_i-\mu_j}.
$$

Let $\la$ be another partition. We write
$\mu\subset\la$ if $\mu_i\le\la_i$ for $i=1,2,\ldots$.

We now prove an important property of the functions $P^*_\la$.
Next, we write $P^*_\mu(\la)$ instead of
$P^*_\mu(\la_1,\ldots,\la_{l(\la)})$.

\proclaim{Theorem 1.5 {\rm (vanishing property)}}{\it

{\rm a)} If $\mu\not\subset\la$, then $P^*_\mu(\la)=0$;

{\rm b)} $P^*_\mu(\mu)=H(\mu)$.
}\endproclaim

(The statement of this Theorem is similar to the vanishing property for
$s^*$-functions in [8] and [6].)

\demo{Proof} Note that $(a\d b)=0$, if
$a,b\in\bZ_+$ and $b>a$.
First, we prove a).

Let $\mu\not\subset\la$, then $\la_k<\mu_k$ for some natural
$k$. Let us choose an arbitrary
$n\ge\max(l(\mu),l(\la))$.

By Definition~1.2,
$$
P^*_\mu(\la)=\frac1{(n-l(\mu))!}\,
\sum\limits_{\om\in S(n)}
F_{\mu|n}(\la_{\om(1)},\ldots,\la_{\om(n)}).
$$
For an arbitrary permutation $\om$ from the group $S(n)$, consider
the corresponding term in the sum
$$
F_{\mu|n}(\la_{\om(1)},\ldots,\la_{\om(n)})=
\prod\limits^{l(\mu)}_{i=1}(\la_{\om(i)}\d\mu_i)
\prod\limits_{\Sb i\le l(\mu)\\
i<j\le n\endSb}
\frac{\la_{\om(i)}+\la_{\om(j)}}{\la_{\om(i)}-\la_{\om(j)}}.
$$
There exists a positive integer $r$ such that $r\le k$ and
$\om(r)\ge k$.
Then we have a chain of inequalities
$$
\la_{\om(r)}\le \la_k<\mu_k\le\mu_r,
$$
therefore,
$$
(\la_{\om(r)}\d\mu_r)=0
$$
and
$$
F_{\mu|n}(\la_{\om(1)},\ldots,\la_{\om(n)})=0.
$$

Since the choice of $\om$ is arbitrary,
$P^*_\mu(\la)=~0$.

b) Arguing as above, we see that
$$
P^*_\mu(\mu)=(\mu_1\d\mu_1)\ldots
(\mu_{l(\mu)}\d\mu_{l(\mu)})
\cdot\prod\limits_{\mu_i>\mu_j}
\frac{\mu_i+\mu_j}{\mu_i-\mu_j}=H(\mu).
$$
Theorem is proved. $\square$
\enddemo

Let $\nu$ be an arbitrary partition. We recall (for details, see~[3])
that the Young diagram of a partition $\nu$ is the set of points
$(i,j)\in\bZ^2$ such that $1\le j\le\la_i,\,1\le i\le l(\la)$.
Let us replace each point
with the unit square with the left upper vertex at this point. We
assume that the first coordinate $i$ (the row index)
increases as one goes downwards, and the second coordinate $j$
(the column index)
increases as one goes from left to right.

For example, if $\nu=(6,5,3,1)$, then
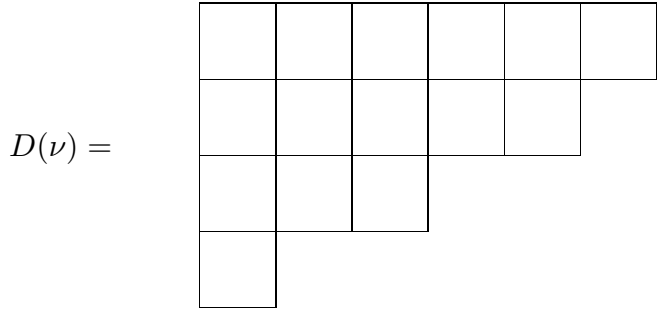
\picture 12 {}
{
\hl 0 0 12
\hl 0 2 12
\hl 0 4 10
\hl 0 6 6
\hl 0 8 2
\vl 0 8 8
\vl 2 8 8
\vl 4 6 6
\vl 6 6 6
\vl 8 4 4
\vl 10 4 4
\vl 12 2 2
\wr -5 4 {$D(\nu)=$}
}
\endpicture

If $\nu$ is a strict partition, then the {\it shifted} diagram
$D'_\nu$ is obtained from the ordinary diagram $D_\nu$ by shifting
the $i$th row $(i-1)$ squares to the right, for all $i>1$. For
$\nu=(6,5,3,1)$, we obtain the shifted diagram
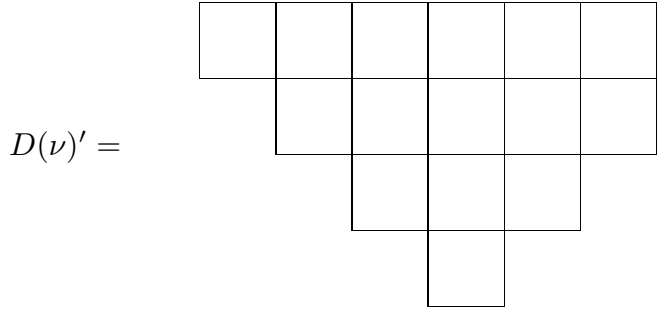
\picture 12 {}
{
\hl 0 0 12
\hl 0 2 12
\hl 2 4 10
\hl 4 6 6
\hl 6 8 2
\vl 0 2 2
\vl 2 4 4
\vl 4 6 6
\vl 6 8 8
\vl 8 8 8
\vl 10 6 6
\vl 12 4 4
\wr -5 4 {$D(\nu)'=$}
}
\endpicture

If $\mu$ and $\la$ are strict partitions, $\mu\subset\la$, then
the {\it skew shifted diagram} $D'_{\la/\mu}$ corresponding to the pair
$(\la,\mu)$ is the difference of the shifted diagrams
$D'_\la$ and $D'_\mu$.

A {\it shifted standard tableau of the form} $\la/\mu$ is a labelling of
the skew shifted diagram $D'_{\la/\mu}$ with the numbers
$1,2,\ldots,|\la|-|\mu|$ such that the numbers strictly increase
from left to right along each row and down each
column.
The {\it dimension} $g_{\la/\mu}$ of a skew shifted diagram
$D'_{\la/\mu}$ is the number of shifted standard tableaux of the form
$\la/\mu$. We put
$g_{\la/\mu}=0$ if $\la$ does not contain $\mu$.
Also put $g_\la=g_{\la/\{\varnothing\}}$,
i.e. the number of shifted standard tableaux
of the form $\la$. There is an explicit formula for $g_\la$ given in
[3, ch.~III, \S8, example 12],
$$
g_\la=\frac{|\la|!}{\la_1!\ldots\la_{l(\la)}!}
\,\prod\limits_{i<j}\frac{\la_i-\la_j}{\la_i+\la_j}.
$$

Now we state the main result of this work
which allows us to obtain an
explicit formula for $g_{\la/\mu}$.

\proclaim{Theorem 1.6}{\it  Suppose
$\mu$ and $\la$ be strict partitions.
Let $m=|\mu|$, $k=|\la|$; then
$$
g_{\la/\mu}=g_\la\cdot\frac{P^*_\mu(\la)}
{(k\d m)}.
$$
}\endproclaim

\demo{Proof} Let us fix $\mu$ and consider
$g_{\la/\mu}$ as a function of strict partition
$\la$, where $|\la|\ge|\mu|$.

For two partitions $\mu$ and $\nu$, we write $\mu\nearrow\nu$
if $|\nu|=|\mu|+1$ and
$\mu\subset\nu$, or, in other words, $D_\nu$ is obtained from
$D_\mu$ by adding one square.

As a function of $\la$, the expression
$g_{\la/\mu}$ is defined by three properties:

(i) $g_{\mu/\mu}=1$;

(ii) if $\la\not\supset\mu$, $|\la|\ge|\mu|$, then
$$
g_{\la/\mu}=0;
$$

(iii)
$g_{\la/\mu}=\sum\limits_{\Sb\nu\text{ is strict,}\\
\nu\nearrow\la\endSb} g_{\nu/\mu}$, $|\la|\ge|\mu|+1$.

Let us prove that these three properties are satisfied for
$$
G(\la)=P^*_\mu(\la)\cdot\frac{g_\la}{(k\d m)}.
$$

(Since $k\ge m$, the denominator does not vanish.)

The property (i) is satisfied,
$$
G(\mu)=P^*_\mu(\mu)\cdot\frac{g_\mu}{m!}=
H(\mu)\cdot\prod\limits_{i<j}
\frac{\mu_i-\mu_j}{\mu_i+\mu_j}
\frac1{\mu_1!\ldots\mu_{l(\mu)}!}=1.
$$

The property (ii) follows from
the vanishing property (Theorem~1.5) for
$P^*_\mu$. Let us prove (iii).
In our case, $k\ge m+1$.

If $\mu\not\subset\la$, then both parts of the desired equality are zero.

Now let $\mu\subset\la$. We set
$$
F(\la_1,\ldots,\la_{l(\la)})=
\frac{(k-m)!}{\la_1!\ldots\la_{l(\la)}!}\,\prod
\limits^{l(\mu)}_{t=1}(\la_t\d\mu_t)
\prod\limits_{l(\mu)<i<j}\frac{\la_i-\la_j}{\la_i+\la_j}.
$$
It follows from Definition~1.2 and the formula for $g_\la$ that
$$
G(\la)=\frac1{(l(\la)-l(\mu))!}\sum\limits_{\om\in S(n)}\sgn(\om)
F(\la_{\om(1)},\ldots,\la_{\om(l(\la))}).
\tag1.3
$$

Let $\la^{(i)}$ denote the partition obtained from
$\la$ by decreasing the $i$th part by 1,
$$
\la^{(i)}=(\la_1,\ldots,\la_{i-1},\la_i-1,\la_{i+1},\ldots,
\la_n).
$$

Let $\om$ be an arbitrary permutation from the group
$S(n)$. The equality
$$
F(\la_{\om(1)},\ldots,\la_{\om(n)})=
\sum\limits^n_{i=1}F(\la^{(i)}_{\om(1)},\ldots,\la^{(i)}_{\om(n)})
$$
is equivalent to the identity
$$
\sum\limits^{l(\la)}_{i=l(\mu)+1}\la_{\om(i)}\!=\!
\sum\limits^{l(\la)}_{i=l(\mu)+1}\la_{\om(i)}\!\!\!
\prod\limits_{\Sb j\ne i,\\
l(\mu)<j\le l(\la)\endSb}\!\!
\frac{\la_{\om(i)}+\la_{\om(j)}}
{\la_{\om(i)}+\la_{\om(j)}-1}\cdot
\frac{\la_{\om(i)}-\la_{\om(j)}-1}
{\la_{\om(i)}-\la_{\om(j)}},
$$
which is given (in another notation for variables) in
[3, ch.~III, \S8, example~12] and [4].

Hence,
we have
for $G(\la)$
$$
G(\la)=\sum\limits_{\nu\nearrow\la} G(\nu).
$$
Note that if $\nu$ is a non--strict partition, then
$$
G(\nu)=0,
$$
because the sum in (1.3) breaks down into pairs of
summands with equal absolute values which occur with different signs.
Therefore, we obtain the relation
$$
G(\la)=\sum\limits_{\Sb\nu\nearrow\la,\\
\nu\text{ is strict}\endSb} G(\nu).
$$

Thus, all three properties are satisfied for $G(\la)$, and Theorem is proved.
$\square$
\enddemo

If $\nu$ and $\eta$ are ordinary partitions (not necessarily strict),
$\nu\subset\eta$, then the skew diagram $D_{\eta/\nu}$ corresponding
to the pair $(\nu,\eta)$ is the difference of diagrams
$D_\eta$ and $D_\nu$.

A standard tableau of the form $\eta/\nu$ is a labelling of squares of the
skew diagram $D_{\eta/\nu}$ with the numbers
$1,2,\ldots,|\eta|-|\nu|$ such that the numbers strictly increase
from left to right along each row and down each
column.
The dimension $f_{\eta/\nu}$ of a skew diagram $D_{\eta/\nu}$ is the
number of standard tableaux of the form $\eta/\nu$.

In [8] there is an explicit formula for $f_{\eta/\nu}$ in terms of so-called
shifted Schur polynomials
$s^*_\mu$. If $l(\mu)\le n$, then
$$
s^*_\mu(x_1,\ldots,x_n)=
\frac{\det[((x_i+n-i)\d(\mu_j+n-j))]^n_{i,j=1}}
{\prod\limits_{i<j}((x_i-x_j)+j-i)}.
$$
(We note that the polynomial $s^*_\mu$, like ordinary Schur polynomials,
possesses the stability property as $n\to\infty$). The formula for
$f_{\eta/\nu}$ takes the form
$$
f_{\eta/\nu}=f_\eta\cdot\frac{s^*_\nu(\eta)}{(|\eta|\d|\nu|)}.
$$
There is an explicit formula
for $f_\eta$
$$
f_\eta=\frac{|\eta|!}
{\prod\limits^{l(\eta)}_{i=1}(\eta_i+l(\eta)-i)!}\cdot
\prod\limits_{1\le i<j\le l(\eta)}
(\eta_i-\eta_j+j-i).
$$

\medskip
\noindent{\bf Remark 1.7.} If $\mu$ and $\la$ are strict partitions,
$\mu\subset\la$, then the skew shifted diagram
$D'_{\la/\mu}$ coincides with an ordinary skew diagram
$D_{\eta/\nu}$ for some partitions $\eta$ and $\nu$ if and only if
one of the following conditions holds:
$$
l(\la)=l(\mu)\quad\text{or}\quad l(\la)=l(\mu)+1.
$$
This means that $D'_\mu$ completely contains the part of the diagram
$D'_\la$ that lies to the left of the vertical line
$j=l(\la)-1$.

If $l(\la)=l(\mu)$ or $l(\la)=l(\mu)+1$, then we define
$\eta$ and $\nu$ as
$$
\aligned
&\eta_i=\la_i+i-1,\\
&\nu_i=\mu_i+i-1,\quad i=1,2,\ldots,l(\la).
\endaligned
$$
Then $D'_{\la/\mu}=D_{\eta/\nu}$.
In this case, the number of standard tableaux
$f_{\eta/\nu}$ equals the number of shifted standard tableaux
$g_{\la/\mu}$.

This implies the following identity for
$P^*_\mu(\la)$ and $s^*_\nu(\eta)$:
$$
P^*_\mu(\la)\cdot\frac{g_\la}{(|\la|\d|\mu|)}=
g_{\la/\mu}=f_{\eta/\nu}=
s^*_\nu(\eta)\cdot
\frac{f_\eta}{(|\eta|\d|\nu))}.
$$
One can check this identity directly from definitions of
$P^*_\mu$ and $s^*_\nu$.

One can compare this identity with a well--known fact from the theory
of Schur superfunctions, the Berele-Regev formula, see
[3, ch.~I, \S3, example~23,(4)].

%2
\head \S2. Proof of the formula for characters of the infinite
spin--symmetric group
\endhead

The symmetric group $S(n)$ (the group of all permutations of the numbers
\linebreak $1,\ldots,n)$
is generated by the permutations $S_k$ of the numbers
$k$ and $k+1$ $(k=1,\ldots,n-1)$ with relations
$$
S^2_k=e;\quad
(S_k\cdot S_{k+1})^3=e;\quad
(S_kS_{k'})^2=e,\quad
k-k'>1.
$$

The group $S(\infty)$ (the group of all finite permutations of
natural numbers) is the inductive limit of the sequence
$S(1)\subset S(2)\subset\ldots$.

The spin--symmetric group $\widetilde S(n)$ is a non--trivial
central $\Bbb Z_2$--extension of the group $S(n)$. It is defined as
the group with generators
$c, t_1, t_2,\ldots,t_{n-1}$ and relations
$$
\gather
c^2=e; \quad ct_k=t_kc;\quad t^2_k=e;\quad(t_kt_{k+1})^3=c;\\
(t_kt_{k'})^2=c,\quad k'-k>1.
\endgather
$$

Projective characters of the group $S(n)$ are linearized by the group
$\widetilde S(n)$, see \cite{4}, \cite{10}, \cite{11}.

Define $\widetilde S(\infty)$ as the inductive limit of the
chain $\widetilde S(1)\subset \widetilde S(2)\subset\ldots$.
The group $\widetilde S(\infty)$ is a non--trivial central
$\Bbb Z_2$--extension of the group $S(\infty)$.

Next, we consider irreducible representations of the group
$\widetilde S(n)$ that send the element $c$  to $-1$.
They may be indexed by strict partitions of the number
$n$ [10]. If a strict partition $\la\vdash n$ is such that
$n-l(\la)$ is even, then there is one irreducible representation
that corresponds to this partition; we denote the character of this
representation by
$\chi^\la$. If a strict partition
$\la\vdash n$ is such that $n-l(\la)$ is odd, then there are two
irreducible representations that correspond to this partition; we denote
their characters by
$\chi^\la_+$ and $\chi^\la_-$.

If $\rho$ is a partition of $n$, then we set
$$
t_\rho=(t_1t_2\ldots t_{\rho_1-1})
(t_{\rho_1+1}\ldots t_{\rho_1+\rho_2-1})\ldots
(t_{\rho_1+\ldots+\rho_{l(\rho)-1}+1}\ldots t_{n-1}).
$$

\noindent
Clearly, it suffices to define characters of irreducible
representations on the elements
$t_\rho$ for their complete definition on the whole group
$\widetilde S(n)$.

In the sequel, $\la,\mu,\nu$ always denote strict partitions.

If $n-l(\la)$ is even, $\la\vdash n$, then it follows from
$\chi^\la(t_\rho)\ne0$ that
$\rho$ is a partition of $n$ into odd parts. If
$n-l(\la)$ is odd, $\la\vdash n$,
then irreducible characters $\chi^\la_+$ and $\chi^\la$ are such that
$\chi^\la_+(t_\rho)$ or $\chi^\la_-(t_\rho)$ can be nonzero only
in the following cases: either all parts of
$\rho$ are odd or
$\rho=\la$. In the first case, $\chi^\la_+(t_\rho)=\chi^\la_-(t_\rho)$,
and, in the second case, $\chi^\la_+(t_\la)=-\chi^\la_-(t_\la)$.

Our purpose is to describe indecomposable characters of
the group $\widetilde S(\infty)$ in the sense of the following definition.

\proclaim{Definition 2.1}{\it  Let $G$ be an abstract group. A function
$\chi:G\to\Bbb C$ is said to be an {\it indecomposable character} if the
following conditions hold:

$1)$ $\chi(e)=1$;

$2)$ $\chi(g_1g_2)=\chi(g_2g_1)$, $\forall g_1$, $g_2\in G$;

$3)$ $\sum\limits_{k,l}\chi(g^{-1}_k g_l)\bar a_k a_l\ge0$,
$\forall g_1,\ldots,g_n\in G$, $\forall a_1,\ldots,
a_n\in\Bbb C$;

$4)$ if $\chi,\chi_1$ and $\chi_2$ satisfy conditions~1)--3), and there
exists a number $a$ such that $0<a<1$ and
$\chi= a\chi_1+(1-a)\chi_2$, then
$\chi=\chi_1=\chi_2$.
}\endproclaim

If $G$ is a finite group, then its indecomposable characters
(in the sense of Definition~2.1) coincide with normalized irreducible
characters.

It follows from Vershik--Kerov theorem [1, 2, ch.~I, \S1] that every
indecomposable character of the group
$\widetilde S(\infty)$ is a pointwise limit of normalized irreducible
charactres of the groups
$\widetilde S(n)$ as $n\to\infty$.

Consider a sequence $(\la(n)\vdash n)$ of strict partitions,
$n=1,2,\ldots$. For every $n$, we denote by
$\xi_n$ any of the following normalized irreducible characters,
$$
\frac{\chi^{\la(n)}}{\chi^{\la(n)}(e)},\quad
\frac{\chi^{\la(n)}_+}{\chi^{\la(n)}_+(e)},\quad
\frac{\chi^{\la(n)}_-}{\chi^{\la(n)}_-(e)}.
\tag2.1
$$

\proclaim{Theorem 2.2 {\rm (M.~L.~Nazarov)}}{\it  The pointwise limit
$\lim\limits_{n\to\infty}\xi_n$ exists if and only if
the limits $\lim\limits_{n\to\infty}\frac{\la_i(n)}n=
\ga_i$, $i=1,2,\ldots$, exist.
}\endproclaim

\demo{Proof} A projection $\pi:\widetilde S(n)\to S(n)$ is defined
as
$$
\pi(c)=e,\quad \pi(t_k)=S_k\quad (k=1,\ldots,n-1).
$$

If at least one cycle in $\pi(t)$ has even length, then
$\xi_n(t)=0$ for all sufficiently large
$n$. It follows that the existence of the limit
$\lim\limits_{n\to\infty}\xi_n$ and its value depend only on the sequence
$\la(n)$, but not on the choice of signs
``$+$'' or ``$-$'' in (2.1).

Thus, we introduce the following notation.
$$
\chi^\la_*=\cases
\chi^\la,\quad &\text{if}\quad |\la|-l(\la)\quad\text{is even},\\
\frac{\chi^\la_++\chi^\la_-}{\sqrt2},\quad&\text{if}\quad
|\la|-l(\la)\quad\text{is odd}.
\endcases
$$

In the space of functions $f(t)$ on the group $\widetilde S(k)$ such that
$f(ct)=-f(t)$, we introduce a scalar product
$$
\l f,g\r_k=\frac1{k!}\sum\limits_{s\in S(k)} (f\bar g)(s).
$$

Denote by $\Res_k$ the operator of restriction to the subgroup
$\widetilde S(k)\subset \widetilde S(n)$, $n>k$.

The pointwise convergence of $\xi_n$ is equivalent
to the following statement.
For every $k$ and every strict partition
$\mu\vdash k$, there exists the limit
$$
\lim\limits_{n\to\infty}\bigg\langle
\frac{\Res_k\chi^{\la(n)}_*}{\chi^{\la(n)}_*(e)},
\chi^\mu_*\bigg\rangle _k.
$$

We denote
$$
\ep(n)=\cases
1,\quad&\text{if}\quad n\quad \text{is odd},\\
\sqrt 2 ,\quad &\text{if}\quad n\quad\text{is even}.
\endcases
$$

The branching rule for $\chi^\nu_*$ takes the following form
[4, Theorem~10.2]
$$
\Res_{|\nu|-1}\chi^\nu_*=\sum\limits_{\mu\nearrow\nu}
\chi^\mu_*\cdot\ep(l(\nu)-l(\mu)),
$$
where we assume that
$\mu$ is a strict partition;
the notation $\mu\nearrow\nu$ was introduced in the proof of the formula
for the dimension of a skew shifted diagram.

Let
$$
\chi^\la_0=2^{\frac{l(\la)-|\la|}2}\,\chi^\la_*=
\cases
2^{\frac{l(\la)-|\la|}2}\,\chi^\la,\,\,&\text{if}\,\,
|\la|-l(\la)\,\,\text{is even},\\
2^{\frac{l(\la)-|\la|-1}2}
(\chi^\la_++\chi^\la_-),\,\,&\text{if}\,\,
|\la|-l(\la)\,\,\text{is odd}.
\endcases
$$

In this normalization, the branching rule takes the simplest form,
$$
\Res_{|\nu|-1}\chi^\nu_0=\sum\limits_{\mu\nearrow\nu}\chi^\mu_0.
$$

This implies the following fact. Suppose $k\le n$, and
$\la\vdash n$ be a strict partition. Then
$$
\Res_k\chi^\la_0=\sum\limits_{\mu\vdash k}g_{\la/\mu}\chi^\mu_0,
$$
where $g_{\la/\mu}$ denotes the dimension of the skew shifted diagram
$D'_{\la/\mu}$.

Taking $\chi^\la_*$ instead of $\chi^\la_0$, we obtain
$$
\Res_k\chi^\la_*=2^{\frac{|\la|-l(\la)}2}\sum\limits_{\mu\vdash k}
g_{\la/\mu}\chi^\mu_0
=\sum\limits_{\mu\vdash k}2^{\frac{|\la|-l(\la)-|\mu|+l(\mu)}2}
g_{\la/\mu}\chi^\mu_*.
$$
Note also that
$$
\chi^\la_*(e)=2^{\frac{|\la|-l(\la)}2}
\chi^\la_0(e)=2^{\frac{|\la|-l(\la)}2}g_\la.
$$

Since $\l\chi^\mu_*,\chi^\mu_*\r_k=1$, we have
$$
\bigg\langle\frac{\Res_k\chi^\la_*}{\chi^\la_*(e)},\chi^\mu_*\bigg\rangle
_k=
\frac{g_{\la/\mu}}{g_\la}\cdot 2^{\frac{l(\mu)-|\mu|}2}
=2^{\frac{l(\mu)-|\mu|}2}\cdot \frac{P^*_\mu(\la)}{(|\la|\d|\mu|)},
$$
where the last equality follows from Theorem~1.6.

Note that Corollary~1.4 implies
$$
P^*_\mu(\la(n))=P_\mu(\la(n))+O(|\la(n)|^{|\mu|-1})
$$
as $n\to\infty$.

Also we have
$$
(n\d|\mu|)=n^{|\mu|}+O(n^{|\mu|-1})
$$
as $n\to\infty$.

Hence, the existence of the pointwise limit
$\lim\limits_{n\to\infty}\xi_n$
is equivalent to the existence of the limits
$$
\lim\limits_{n\to\infty}\frac{P_\mu(\la(n))}
{n^{|\mu|}}
\tag2.2
$$
for all $\mu$.

The functions $P_\mu$ form a linear basis of the algebra of
supersymmetric functions $\Om$. Since odd Newton sums
$p_1,p_3,p_5,\ldots$ generate algebraically
$\Om$, the existence of the limits (2.2) is equivalent to the existence
of the limits
$$
\lim\limits_{n\to\infty}\frac{p_m(\la(n))}{n^m},\quad
m=1,3,5,\ldots.
$$

We will prove that if there exist the limits
$$
\lim\limits_{n\to\infty}\frac{\la_i(n)}n=\ga_i\quad
(i=1,2,\ldots),
$$
then there exist the limits
$$
\lim\limits_{n\to\infty}\frac{p_m(\la(n)}{n^m}=
\sum\limits^\infty_{k=1}\ga^m_k=
p_m(\ga_1,\ga_2,\ldots)
$$
for $m=3,5,7,\ldots$. $\left(\text{Note that} \,\,
\frac{p_1(\la(n))}n=1\right)$ .

If $l(\la(n))$ is bounded as $n\to\infty$, then this statement is trivial.

Hence, we may assume that $\lim\limits_{n\to\infty}
l(\la(n))=\infty$. Then the desired fact is evident from the
following estimation. Suppose
$\Cal N\le l(\la)$, $m\ge 3$. Then
$$
\multline
\frac{p_m(\la_1,\ldots,\la_{l(\la)})}{n^m}=
\frac{p_m(\la_1,\ldots,\la_{\Cal N})}{n^m}+
\sum\limits^{l(\la)}_{i=\Cal N+1}\frac{\la^m_i}{n^m}\le\\
\le
\frac{p_m(\la_1,\ldots,\la_{\Cal N})}{n^m}+\frac{\la_{\Cal N+1}}
n.
\endmultline
$$
In the last inequality we used the fact that
$$
\sum\limits^{l(\la)}_{i=\Cal N+1}\la^{m-1}_i\le
\sum\limits^{l(\la)}_{i=1}\la^{m-1}_i\le
\left(\sum\limits^{l(\la)}_{i=1}\la_i\right)^{m-1}=
n^{m-1}.
$$
On the other hand,
$$
\frac{p_m(\la_1,\ldots,\la_{l(\la)})}{n^m}\ge
\frac{p_m(\la_1,\ldots,\la_{\Cal N})}{n^m}.
$$
Taking into account that
$$
\lim\limits_{r\to\infty}\,\sum\limits^\infty_{k=r+1}
\ga^m_k=0,\quad m=3,5,7,\ldots,
$$
and
$$
\lim\limits_{\Cal N\to\infty}\left(
\lim\limits_{n\to\infty}\frac{\la_{\Cal N+1}(n)}
n\right)=\lim\limits_{\Cal N\to\infty}\ga_{\Cal N+1}=0,
$$
we obtain the desired statement.

Thus, the existence of the limits
$$
\lim\limits_{n\to\infty}
\frac{\la_i(n)}n\quad(i=1,2,\ldots)
$$
implies the existence of the pointwise limit
$\lim\limits_{n\to\infty}\xi_n$.

Now we will prove the inverse statement. We still consider
a sequence of strict partitions
$(\la(n))_{n\ge1}$ such that
$|\la(n)|=n$.

As proved above,
the existence of the pointwise limit
$\lim\limits_{n\to\infty}\xi_n$ implies the existence of the limits
$$
\lim\limits_{n\to\infty}\frac{p_m(\la(n))}{n^m},\quad
m=1,3,5,\ldots.
$$
The sequences $\left(\frac{\la_i(n)}n\right)_{n\ge1}$ are bounded for all
$i$. We can choose a subsequence
$(\la(n_k))_{k\ge1}$ such that the limits
$$
\lim\limits_{k\to\infty}\frac{\la_i(n_k)}{n_k}
$$
exist for all $i$. Let
$$
\lim\limits_{k\to\infty}\frac{\la_i(n_k)}{n_k}=\ga_i,\quad
i=1,2,\ldots.
$$
Then, arguing as above, we obtain the equalities
$$
\lim\limits_{n\to\infty}\frac{p_m(\la(n))}{n^m}=
\lim\limits_{k\to\infty}\frac{p_m(\la(n_k))}{n^m_k}=
p_m(\ga_1,\ga_2,\ldots)
$$
for $m=3,5,7,\ldots$.
Note that every sequence
$\ga=(\ga_i)_{i\ge1}$ such that $\ga_1\ge\ga_2\ge\ldots
\ge0$, $\sum\limits_i\ga_i\le 1$
is uniquely defined by the values
$$
p_3(\ga),p_5(\ga),\ldots\,.
$$
In fact, $\ga_1$ is uniquely defined by the condition
that
$$
\lim\limits_{n\to\infty}\frac{p_{2n+1}(\ga)}{\ga_1^{2n+1}}
$$
exists and is not zero. In order to define
$\ga_{i+1}$, given $\ga_1,\ga_2,\ldots,\ga_i$,
it suffices to consider the values
$$
p_3(\ga)-\sum\limits^i_{k=1}\ga^3_k,
p_5(\ga)-\sum\limits^i_{k=1}\ga^5_k,\ldots\,.
$$
It follows that for all $i$, the sequence
$\left(\!\frac{\la_i(n)}n\!\right)_{n\ge1}$ has only one limit point,
and thus it converges. Theorem follows.
$\square$
\enddemo

Now we will find the pointwise limit of the sequence
$\xi_n$ as $\lim\limits_{n\to\infty}\frac{\la_i(n)}n=\ga_i$,
$i=1,2,\ldots$.

\proclaim{Theorem 2.3 {\rm (M.~L.~Nazarov [5])}}{\it
Suppose
$$
\ga=(\ga_1\ge\ga_2\ge\ldots\ge0)
$$
and
$$
\lim\limits_{n\to\infty}\frac{\la_i(n)}n=\ga_i,
\quad i=1,2,\ldots,\la(n)\vdash n.
$$
We denote by $\psi_\ga$ the pointwise limit of the sequence
$\xi_n$ defined by $(2.1)$. Then
$$
\psi_\ga(t_\rho)=
\cases
\prod\limits_{i\ge2}p_i(\ga)^{m_i(\rho)}\cdot
2^{\frac{l(\rho)-|\rho|}2},\quad&\text{if}\quad\rho
\,\,\text{is a partition}\\
&\text{ into odd parts},\\
0,\quad&\text{otherwise},
\endcases
$$
where $m_i(\rho)$ is the number of parts of the partition $\rho$
equal to $i$, and $t_\rho$ is the element of the group
$\widetilde S(\infty)$ that was introduced above.
}\endproclaim

\demo{Proof} Let $\rho$ be a partition of $k$ into
odd parts. Then
$$
\psi_\ga(t_\rho)=
\lim\limits_{n\to\infty}
\frac{\chi^{\la(n)}_*(t_\rho)}
{\chi^{\la(n)}_*(e)}=
\lim\limits_{n\to\infty}\sum\limits_{\mu\vdash k}
\left\l\frac{\Res_k\chi^{\la(n)}_*}
{\chi^\la_*(e)},\,
\chi^\mu_*\right\r
_k\chi^\mu_*(t_\rho).
\tag2.3
$$
Taking into account the above--obtained equality
$$
\left\l\frac{\Res_k\chi^\la_*}{\chi^\la_*(e)},\chi^\mu_*\right\r_k=
\frac{2^{\frac{l(\mu)-k}2}P^*_\mu(\la)}{(|\la|\d k)},
$$
we rewrite
(2.3) as
$$
%\multline
\lim\limits_{n\to\infty}\sum\limits_{\mu\vdash k}
\frac{2^{\frac{l(\mu)-k}2}P^*_\mu(\la(n))}
{(|\la(n)|\d k)}\cdot \chi^\mu_*(t_\rho)=
\lim\limits_{n\to\infty}
\sum\limits_{\mu\vdash k}
\frac{2^{\frac{l(\mu)-k}2}\cdot P_\mu(\la(n))}
{|\la(n)|^k}\cdot
\chi^\mu_*(t_\rho)
%\endmultline
\tag2.4
$$
Now we use an equality from
[11,~\S7] which in our notation takes the form
$$
\sum\limits_{\mu\vdash k}
2^{\frac{l(\mu)-k}2}
P_\mu(\la)\cdot \chi^\mu_*(t_\rho)
=2^{\frac{l(\rho)-k}2}p_\rho(\la).
$$
Hence (2.4) takes a simple form
$$
\lim\limits_{n\to\infty}\frac{p_\rho(\la(n))}
{|\la(n)|^k}\cdot 2^{\frac{l(\rho)-k}2}=
2^{\frac{l(\rho)-|\rho|}2}
\prod\limits_{i\ge2}
p_i(\ga)^{m_i(\rho)}.
$$

If at least one part of the partition $\rho$ is even, then
$$
\psi_\ga(t_\rho)=0,
$$
since then
$$
\xi_n(t_\rho)=0
$$
for all sufficiently large $n$.
This completes the proof.
$\square$
\enddemo

\proclaim{Proposition 2.4}{\it  The functions $\psi_\ga$
obtained in the previous Theorem are indeed indecomposable
characters of the group
$\widetilde S(\infty)$.
}\endproclaim

\demo{Proof}
We check Definition~2.1 for the functions
$\psi_\ga$. Properties 1)--3) are satisfied since they are satisfied
for normalized characters whose limit is
$\psi_\ga$. For $\psi_\ga$, the multiplicativity property
is satisfied, i.e., if
$\rho$ and $\si$ are arbitrary partitions,
and $\rho\cup\si$ denotes their disjoint union, then
$$
\psi_\ga(t_{\rho\cup\si})=\psi_\ga(t_\rho)\cdot\psi_\ga
(t_\si).
$$
This implies the indecomposability of
$\psi_\ga$
(property~4 of Definition~2.1) [13, 7].
\enddemo
Supported by Soros International Educational Program,
grant 2093s.

Translated by N.~V.~Tsilevich.

\enddocument